\documentclass[psamsfonts]{amsart}  % proc-l
\usepackage{amssymb}
\usepackage[dvips]{graphicx}

\usepackage{graphicx}
\usepackage{amssymb}
\usepackage{amsmath}
\usepackage{color}
\usepackage{times}

\usepackage[unicode,bookmarks,colorlinks]{hyperref}
\hypersetup{
    linkcolor=brickred,
}

\definecolor{mahogany}{cmyk}{0, 0.77, 0.87, 0}
\definecolor{salmon}{cmyk}{0, 0.53, 0.38, 0}
\definecolor{melon}{cmyk}{0, 0.46, 0.50, 0}
\definecolor{yellowgreen}{cmyk}{0.44, 0, 0.74, 0}
\definecolor{brickred}{cmyk}{0, 0.89, 0.94, 0.28}
\definecolor{OliveGreen}{cmyk}{0.64, 0, 0.95, 0.40}
\definecolor{RawSienna}{cmyk}{0, 0.72, 1.0, 0.45}
\definecolor{ZurichRed}{rgb}{1, 0, 0} % Red of svgnames

\usepackage{fancyhdr}
\usepackage{amsmath,amstext,amssymb,amsopn,amsthm}
\usepackage{amsmath,amssymb,amsthm}
\usepackage[mathscr]{eucal}

\begin{document}

%\newtheorem{thm}{Theorem}
%\numberwithin{thm}{Theorem}
\newtheorem{lemma}[thm]{Lemma}
\newtheorem{remark}{Remark}
%\newtheorem{corr}[thm]{Corollary}[sections]
%\numberwithin{corollary}[thm]{section}
\newtheorem{proposition}{Proposition}
\newtheorem{theorem}{Theorem}[section]
\newtheorem{deff}[thm]{Definition}
\newtheorem{case}[thm]{Case}
%\numberwithin{deff}{section}
\newtheorem{prop}[thm]{Proposition}
%\numberwithin{equation}{subsection}
%\numberwithin{equation}{section}
\newtheorem{example}{Example}

\newtheorem{corollary}{Corollary}

\numberwithin{equation}{section}
\numberwithin{definition}{section}
%\numberwithin{problem}{section}
\numberwithin{corollary}{section}
%\numberwithin{proposition}{subsection}

\numberwithin{theorem}{section}

\numberwithin{remark}{section}
\numberwithin{example}{section}
\numberwithin{proposition}{section}

\newcommand{\gap}{\lambda_{2,D}^V-\lambda_{1,D}^V}
\newcommand{\gapR}{\lambda_{2,R}-\lambda_{1,R}}
\newcommand{\bD}{\mathrm{I\! D\!}}
\newcommand{\calD}{\mathcal{D}}
\newcommand{\calA}{\mathcal{A}}

\newcommand{\conjugate}[1]{\overline{#1}}
\newcommand{\abs}[1]{\left| #1 \right|}
\newcommand{\cl}[1]{\overline{#1}}
\newcommand{\expr}[1]{\left( #1 \right)}
\newcommand{\set}[1]{\left\{ #1 \right\}}

\newcommand{\calC}{\mathcal{C}}
\newcommand{\calE}{\mathcal{E}}
\newcommand{\calF}{\mathcal{F}}
\newcommand{\Rd}{\mathbb{R}^d}
\newcommand{\BR}{\mathcal{B}(\Rd)}
\newcommand{\R}{\mathbb{R}}
\newcommand{\al}{\alpha}
\newcommand{\RR}[1]{\mathbb{#1}}
\newcommand{\bR}{\mathrm{I\! R\!}}
\newcommand{\ga}{\gamma}
\newcommand{\om}{\omega}
\newcommand{\A}{\mathbb{A}}
\newcommand{\bH}{\mathbb{H}}

\newcommand{\bb}[1]{\mathbb{#1}}
\newcommand{\bI}{\bb{I}}
\newcommand{\bN}{\bb{N}}

\newcommand{\uS}{\mathbb{S}}
\newcommand{\M}{{\mathcal{M}}}
\newcommand{\calB}{{\mathcal{B}}}

\newcommand{\W}{{\mathcal{W}}}

\newcommand{\m}{{\mathcal{m}}}

\newcommand {\mac}[1] { \mathbb{#1} }

\newcommand{\bC}{\Bbb C}

\newtheorem{rem}[theorem]{Remark}
\newtheorem{dfn}[theorem]{Definition}
\theoremstyle{definition}
\newtheorem{ex}[theorem]{Example}
\numberwithin{equation}{section}

\newcommand{\Pro}{\mathbb{P}}
\newcommand\F{\mathcal{F}}
\newcommand\E{\mathbb{E}}
\newcommand\e{\varepsilon}
\def\H{\mathcal{H}}
\def\t{\tau}

\newcommand{\blankbox}[2]{%
  \parbox{\columnwidth}{\centering
%    Set fboxsep to 0 so that the actual size of the box will match the
%    given measurements more closely.
    \setlength{\fboxsep}{0pt}%
    \fbox{\raisebox{0pt}[#2]{\hspace{#1}}}%
  }%
}
%\copyrightinfo{2007}{American Mathematical Society}
%\begin{document}

%\copyrightinfo{2007}{American Mathematical Society}
%\begin{document}
%\maketitle
\title[Fourier multipliers]{Martingales and Sharp Bounds for Fourier multipliers}\thanks{To appear in Annales Academiae Scientiarum Fennicae Mathematica}

\author{Rodrigo Ba\~nuelos}\thanks{R. Ba\~nuelos is supported in part  by NSF Grant
\# 0603701-DMS}
\address{Department of Mathematics, Purdue University, West Lafayette, IN 47907, USA}
\email{banuelos@math.purdue.edu}
\author{Adam Os\c ekowski}\thanks{A. Os\c ekowski is supported in part by MNiSW Grant N N201 364436}
\address{Department of Mathematics, Informatics and Mechanics, University of Warsaw, Banacha 2, 02-097 Warsaw, Poland}
\email{ados@mimuw.edu.pl}

\begin{abstract} 
Using the argument of  Geiss, Montgomery-Smith and Saksman \cite{GMSS}, and a new martingale  inequality, the $L^p$--norms of certain Fourier multipliers in $\R^d$, $d\geq 2$, are identified.   These include, among others,  the second order Riesz transforms $R_j^2$, $j=1, 2, \dots, d$, and some of the L\'evy multipliers studied in \cite{BBB}, \cite{BB}. 

\end{abstract}
\maketitle

\section{Introduction}
Martingale inequalities have played a fundamental role for many years in obtaining bounds for the $L^p$-norms of many important singular integrals and Fourier multipliers, both in the real setting and in the Banach space setting.  At the root of these results are the fundamental inequalities of Burkholder on martingale transforms. There is now a huge literature on this subject which would be impossible to review here.   For an overview of this literature, see \cite{Ban} and \cite{BanDav}. The purpose of this paper is to show that there are several instances where some of the  upper bounds, and especially those obtained in recent years, are also lower bounds, hereby enlarging the class of Fourier multipliers where one can compute the norms exactly.   These results are motivated by the paper of Geiss, Montgomery-Smith and Saksman \cite{GMSS}, which has its roots in the work of Bourgain \cite{Bou1}.  The Bourgain result itself is also rooted in the inequalities of Burkholder. While our proof of Theorem \ref{main} is a small modification of  the Geiss, Montgomery-Smith, Saksman argument, we believe our results here will further stimulate interest on these problems and their connections to the (still open) celebrated conjecture of Iwaniec \cite{Iwa} concerning the norm of the Beurling-Ahlfors operator. See \cite{Ban} for some of the history and recent results related to this conjecture.

 Let  $f=\{f_n, n\geq 0\}$ be a martingale on a probability space $(\Omega, \mathcal{F}, \mathbb{P})$ 
with respect to the  sequence of $\sigma$-fields $\calF_n\subset \calF_{n+1}$, $n\geq 0$, contained in $\calF$. The sequence  
$df=\{df_k, k\geq 0\}$, where $df_k=f_k-f_{k-1}$ for $k\geq 1$ and $df_0=f_0$, is called the martingale difference sequence of $f$. 
Thus $f_n=\sum_{k=0}^n df_k$ for all $n\geq 0$.  Given a sequence of random variables $\{v_k, k\geq 0\}$ uniformly bounded by $1$ for all $k$ and with each $v_k$  measurable with respect to $\calF_{(k-1)\vee 0}$ (such sequence is said to be predictable), the martingale difference sequence $\{v_kdf_k, k\geq 0\}$ generates a new martingale called the {\it ``martingale transform"} of $f$ and  denoted by $v\ast f$.  Thus $(v\ast f)_n=\sum_{k=0}^n v_kdf_k$ for all $n\geq 0$. The maximal function of a martingale is denoted by $f^*=\sup_{n\geq 0}|f_n|$.  We also set $\|f\|_p=\sup_{n\geq 0}\|f_n\|_p$ for $0<p<\infty$. Burkholder's 1966 result in \cite{Bur} asserts that the operator $f\to v\ast f=g$ is bounded on $L^p$ for all $1<p<\infty$. In his 1984 seminal paper \cite{B0} Burkholder determined the norm of this operator. For $1<p<\infty$ we let $p^*$  denote the maximum of $p$ and $q$, where $\frac{1}{p}+\frac{1}{q}=1$.   Thus $p^* = \max\{p,\frac{p}{p-1}\}$ and 
\begin{equation}\label{p^*}
p^*-1=\begin{cases} \frac{1}{p-1},  \hskip4mm  1<p\leq 2,\\
p-1 , \hskip3mm  2\leq p <\infty.
\end{cases}
\end{equation}
%This notation will be used throughout this paper.  
\begin{theorem}\label{Burkholder1} Let $f=\{f_n, n\geq 0\}$ be a martingale with difference sequence $df=\{df_k, k\geq 0\}$.   Let $g=v\ast f$ be the martingale transform of $f$ by the real predictable sequence $v=\{v_k, k\geq 0\}$ uniformly bounded in absolute value by 1.  Then 
\begin{equation}\label{bur2}
\|g\|_p\leq (p^*-1)\|f\|_p, \qquad  1<p<\infty,
\end{equation}
and the constant $p^*-1$ is best possible. 
\end{theorem}

By considering dyadic martingales, inequality \eqref{bur2} contains the classical inequality of  Marcinkiewicz \cite{Marc} and Paley \cite{Pal} for Paley-Walsh martingales with the optimal constant. 

\begin{corollary}Let $\{h_k,  k\geq 0\}$ be the Haar system in the Lebesgue unit interval $[0, 1)$.  That is,  $h_0=[0, 1), h_1=[0, \,1/2)-[1/2, \,1),  h_3=[0, \,1/4)-[1/4, \,1/2), h_4=[1/2, \,3/4)-(3/4,\, 1), \dots$, where the same notation is used for an interval as for its indicator function. Then for any sequence $\{a_k, {k\geq 0}\}$ of real numbers and any sequence $\{\e_k, {k\geq 0}\}$ of signs,
\begin{equation}\label{Paleyreal}
\Big\|\sum_{k=0}^{\infty} \varepsilon_k a_k h_k\Big\|_p\leq (p^*-1)\Big\|\sum_{k=0}^{\infty}  a_k h_k\Big\|_p, \qquad 1<p<\infty.
\end{equation}
The constant $p^*-1$ is best possible.
\end{corollary}

In \cite{C}, K.P. Choi used the techniques of Burkholder to identify the best constant in the martingale transforms where the predictable sequence $v$ takes values in $[0, 1]$ instead of $[-1, 1]$. While Choi's  constant is not as explicit as the $p^*-1$ constant of Burkholder, one does have a lot of information about it. 

 \begin{theorem}\label{Choi1} Let $f=\{f_n, n\geq 0\}$ be a real-valued martingale with difference sequence $df=\{df_k, k\geq 0\}$.   Let $g=v\ast f$ be the martingale transform of $f$ by a predictable sequence $v=\{v_k, k\geq 0\}$ with values in $[0,1]$.  Then 
\begin{equation}\label{choi1}
\|g\|_p\leq c_p\|f\|_p, \qquad  1<p<\infty,
\end{equation}
with the best constant $c_p$ satisfying
$$
c_p=\frac{p}{2}+ \frac{1}{2}\log\left(\frac{1+e^{-2}}{2}\right) +\frac{\alpha_2}{p}+\cdots
$$
where 
$$\alpha_2=\left[\log\left(\frac{1+e^{-2}}{2}\right)\right]^2+\frac{1}{2}\log\left(\frac{1+e^{-2}}{2}\right)-2\left(\frac{e^{-2}}{1+e^{-2}}\right)^{2}. $$
\end{theorem}
 
As observed by Choi, 
\begin{equation}
c_p\approx \frac{p}{2}+\frac{1}{2}\log\left(\frac{1+e^{-2}}{2}\right),
\end{equation}
with this approximation becoming better for large $p$. 
It also follows trivially from Burkholder's inequalities that  (even without knowing explicitly the best constant $c_p$) 
\begin{equation}\label{choibound}
\max\left(1, \frac{p^*-1}{2}\right)\leq c_p\leq \frac{p^*}{2}.
\end{equation}
As in the case of Burkholder, Choi's result gives

\begin{corollary}\label{Choi2} Let $\{h_k, k\geq 0\}$ be the Haar system as above. Then for any sequence $\{a_k, {k\geq 0}\}$ of real numbers and any sequence $\{\e_k, {k\geq 0}\}$ of numbers in $\{0,1\}$,
\begin{equation}\label{choi2}
\Big\|\sum_{k=0}^{\infty} \varepsilon_k a_k h_k\Big\|_p \leq c_p\Big\|\sum_{k=0}^{\infty}  a_k h_k\Big\|_p, \qquad  1<p<\infty,
\end{equation}
where $c_p$ is the constant in \eqref{choi1}.  The inequality is sharp. 
\end{corollary}

Motivated by Theorems \ref{Burkholder1} and \ref{Choi1} we introduce a new constant. 

\begin{dfn}\label{defC}
Let $-\infty<b<B<\infty$ and $1<p<\infty$ be given and fixed. We define $C_{p,b,B}$ as the least positive number $C$ such that for any real-valued martingale $f$ and for any transform $g=v\ast f$ of $f$ by a predictable sequence $v=\{v_k, k\geq 0\}$ with values in $[b,B]$,  we have
\begin{equation}\label{martin}
 ||g||_p\leq C||f||_p.
\end{equation}
\end{dfn}

Thus, for example, $C_{p,-a,a}=a(p^*-1)$ by Burkholder's Theorem 
 \ref{Burkholder1} and $C_{p,0,a}=a\,c_p$ by Choi's Theorem \ref{Choi1}.  It is also the case that for any $b, B$ as above, 
 $C_{p,b,B}\leq \max\{B, |b|\}(p^*-1)$ and in fact a simple transformation gives that 
 
 \begin{equation}
 \max\left\{\left(\frac{B-b}{2}\right)(p^*-1), \, \max\{|B|, |b|\}\right\} \leq C_{p,b,B}\leq \frac{(B-b)}{2} p^* +b.
 \end{equation}

 We also point out that by a result of Maurey \cite{Mau}, and independently of Burkholder \cite{B-1}, the constant $C_{p,b,B}$ in this definition remains the same if we consider Paley-Walsh martingales only. Furthermore, the reasoning presented in the Appendix of \cite{Bu} shows that if the transforming sequence is deterministic and takes values in $\{b,B\}$, then the constant in \eqref{martin} does not change either.    

A bounded, complex valued function $m$ on $\R^d\setminus\{0\}$, $d\geq 1$, is called a Fourier multiplier. We define the operator $T_m:L^2(\R^d)\to L^2(\R^d)$ associated to $m$ by 
\begin{equation}\label{multipler}
T_mf=\mathcal{F}^{-1}(m\mathcal{F}),
\end{equation}
where $\mathcal{F}$ is a Fourier transform 
$$\mathcal{F}f(\xi)=\widehat{f}(\xi)=\int_{\R^d} e^{-i\langle \xi, x\rangle}f(x)\mbox{d}x.$$
The multiplier $m$ is said to be homogeneous of order $0$ if  $m(\lambda \xi)=m(\xi)$ for all $\xi\in \R^d\setminus \{0\}$ and $\lambda>0$, and it is said to be even if  $m(\xi)=m(-\xi)$ for all $\xi\in \R^d\setminus \{0\}$.  
We will be particularly interested in those $m$ for which the corresponding $T_m$ is bounded on $L^p(\R^d)$, $1<p<\infty$ (more formally, has a bounded extension to $L^p(\R^d)$). To shorten the notation, we will usually denote the operator norm  $||T_m:L^p(\R^d)\to L^p(\R^d)||$ just by $\|T_m\|_p$, when no danger of confusion exists.  

As an application of the above martingale inequalities to Fourier multipliers, we have the following theorem.  

\begin{theorem}\label{main}
Let $d\geq 2$ be a given integer. Let $m$ be a real and even multiplier which is homogeneous of order 0 on $\R^d$. 
 Denote by $b$ and $B$ the minimal and the maximal term of the sequence
 $$ \big(m(1,0,0,\ldots,0),m(0,1,0,\ldots,0),\ldots,m(0,0,\ldots,0,1)\big),$$
respectively. Then for $1<p<\infty$ and $C_{p,b,B}$ as in Definition \ref{defC}, we have
\begin{equation}\label{mainin}
\|T_m\|_p\geq C_{p,b,B}.
\end{equation}
Furthermore,  since $\|T_m\|_p$ is preserved under rotations and reflections of the multiplier, we have 
\begin{equation}\label{mainin2}
\|T_m\|\geq \sup_{e}C_{p,b(e),B(e)}, \qquad 1<p<\infty,
\end{equation}
where the supremum runs over all orthonormal bases $e=(e_j)_{j=1}^d$ of $\R^d$ and $b(e)$, $B(e)$ stand for the minimal and the maximal term of the sequence $m(e_1)$, $m(e_2)$, $\ldots$, $m(e_d)$, respectively.
\end{theorem}

Recall that the Riesz transforms  $R_1$, $R_2$, $\ldots$, $R_d$ in $\bR^d$, $d\geq 2$, are the Fourier multipliers given by
$$ \widehat{R_jf}(\xi)=-i\frac{\xi_j}{|\xi|} \widehat{f}(\xi),\quad \xi\in \R^d\setminus\{0\},\qquad j=1,\,2,\,\ldots,d.$$
These multipliers do not satisfy the assumptions of the above theorem: they are neither real nor even. However, they give rise to the second order Riesz transforms,
$$ \widehat{R_jR_kf}(\xi)=\frac{-\xi_j \xi_k}{|\xi|^2} \widehat{f}(\xi),\quad \xi\in \R^d\setminus\{0\},\quad j,\, k=1,\,2,\,\ldots,d,$$
which have the desired properties. 
It was proved by Nazarov and Volberg \cite{NazVol} and  Ba\~nuelos and M\'endez-Hern\'andez \cite{BMH} that 
\begin{equation}\label{rieszbound1}
||R_1^2-R_2^2||_p=||2R_1R_2||_p\leq C_{p,-1,1}\leq p^*-1.
\end{equation}
Geiss, Montgomery-Smith and Saksman \cite{GMSS} showed that the inequality in the reverse direction is also true and hence
\begin{equation}\label{rieszbound2}
||R_1^2-R_2^2||_p=||2R_1R_2||_p=C_{p,-1,1}=p^*-1.
\end{equation}
We shall establish the following extension of this result.   

\begin{theorem}\label{Riesz1}
Let $d\geq 2$ and assume that $\mathbb{A}=(a_{ij})_{i,j=1}^d$ is a $d\times d$ symmetric matrix  with real entries and eigenvalues $\lambda_1\leq \lambda_2\leq \ldots \leq \lambda_d$.  Consider the operator $\mathcal{S}_{\mathbb{A}}=\sum_{i,j=1}^d a_{ij}R_iR_j$ with the multiplier
$ m(\xi)=\frac{(\mathbb{A}\xi,\xi)}{|\xi|^2}.$  Then for $1<p<\infty$,
\begin{equation}\label{generalriesz}
%\|\sum_{i,j=1}^d a_{ij}R_iR_j\|_p=C_{p,\lambda_1,\lambda_d}.
\|\mathcal{S}_{\mathbb{A}}\|_p=C_{p,\lambda_1,\lambda_d}.
\end{equation}
%This follows immediately from \eqref{mainin2}, applied to the basis $e$.
\end{theorem}

\begin{corollary}\label{Ries2} If $d\geq 2$ and $J\subsetneq \{1,2,\,\ldots,d\}$, then 

\begin{equation}\label{squareriesz}
||\sum_{j\in J} R_j^2||_p= C_{p,0,1}=c_p, \quad 1<p<\infty,
\end{equation}
where $c_p$ is the Choi constant in \eqref{choi1}. 
\end{corollary}

The lower bound in  \eqref{generalriesz} follows from \eqref{mainin2} applied to the basis of eigenvectors 
$(e_1,\,e_2,\,\ldots,\, e_d)$ corresponding to  $\lambda_1\leq \lambda_2\leq \ldots \leq \lambda_d$.  The upper bound follows from the stochastic integral representation for these  operators first introduced in Ba\~nuelos and M\'endez-Hern\'andez \cite{BMH} and the Burkholder--type inequality \eqref{maininmart} below for continuous time martingales under a more general (not necessarily symmetric) subordination condition.  This result is of independent interest and can be applied to the L\'evy multipliers studied in \cite{BBB} and \cite{BB}, as we shall see momentarily. 

To introduce the necessary notions in the continuous-time setting, suppose that $(\Omega,\mathcal{F},\mathbb{P})$ is a complete probability space, filtered by $(\mathcal{F}_t)_{t\geq 0}$, a nondecreasing and right-continuous family of sub-$\sigma$-fields of $\mathcal{F}$. Assume, as usual, that $\F_0$ contains all the events of probability $0$. Let $X$, $Y$ be adapted, real valued martingales which have right-continuous paths with left-limits (r.c.l.l.).  Denote by $[X,X]$ the quadratic variation process of $X$: we refer the reader to Dellacherie and Meyer \cite{DM} for details. 
Following  Ba\~nuelos and Wang \cite{BW} and  Wang \cite{W}, we say that $Y$ is differentially subordinate to $X$ if the process $([X,X]_t-[Y,Y]_t)_{t\geq 0}$ is nondecreasing and nonnegative as a function of $t$. We have the following extension of Theorem \ref{Burkholder1}, proved by Ba\~nuelos and Wang \cite{BW} for continuous-path martingales and by Wang \cite{W} in the general case. Namely, if $Y$ is differentially subordinate to $X$, then
\begin{equation}\label{symmetric}
\|Y\|_p\leq (p^*-1)\|X\|_p, \qquad 1<p<\infty,
\end{equation}
and the inequality is sharp. Here $\|X\|_p$, the $p$-th moment of $X$, is defined analogously as in the discrete time: $\|X\|_p=\sup_{t\geq 0}\|X_t\|$, $0<p<\infty$. The following theorem extends this result and can be regarded as a continuous-time version of the inequality for non-symmetric martingale transforms.
 
\begin{theorem}\label{Non-Symm} Let $-\infty<b<B<\infty$ and suppose that $X$, $Y$ are real valued  martingales satisfying the non--symmetric subordination condition 

\begin{equation}\label{nonsymm}
 \mbox{d}\left[Y-\frac{b+B}{2}X,Y-\frac{b+B}{2}X\right]_t \leq \mbox{d}\left[\frac{B-b}{2}X,\frac{B-b}{2}X\right]_t,
\end{equation}
for all $t\geq 0$. Then 
\begin{equation}\label{maininmart}
 ||Y||_p\leq C_{p,b,B}||X||_p, \quad 1<p<\infty,
\end{equation}
and the inequality is sharp.
\end{theorem}
Let us clarify that for $t=0$, the condition \eqref{nonsymm} means that
$$ \left(Y_0-\frac{B+b}{2}X_0\right)^2\leq \left(\frac{B-b}{2}X_0\right)^2,$$
or $(Y_0-bX_0)(Y_0-BX_0)\leq 0$.  

%We will outline the proof of the upper bound for Theorem \ref{Riesz1} in \S3.  
%It follows immediately from follows %immediately from 
%Theorem \ref{Non-Symm} and the now standard arguments from \cite{BMH}. 
Theorem \ref{Non-Symm} combined with the techniques from \cite{BBB} and 
\cite{BB} yields new results for multipliers arising from L\'evy processes.  
%We leave this now standard argument to the reader.  
%For more applications, see also \cite{Ban}. 
Consider a measure  $\nu\geq 0$ 
on $\bR^d$ satisfying $ \nu(\{0\})=0$ and
 \begin{equation}\label{levymu1}
 \int_{\bR^d} \frac{|x|^2}{1+|x|^2} \,\mbox{d}\nu(x)<\infty.
\end{equation}
A measure with these properties is called  a  L\'evy measure.  For any finite Borel measure  $\mu\geq
0$ on the unit sphere
$\uS\subset \bR^d$ and any functions $\varphi:\bR^d\to \bC$, $\psi:\uS\to \bC$ with $\|\phi\|_\infty\leq 1$ and $\|\psi\|_\infty\leq 1$, 
we consider the multiplier

\begin{equation} \label{defmuL}
m\left(\xi\right)=
\frac{
      \int_{\bR^d} \Big(1 - \cos \langle \xi, x\rangle \Big)\varphi\left(x\right)\mbox{d}\nu(x)+ 
     \frac{1}{2} \int_{\uS} |\langle\xi, \theta\rangle|^2 \psi\left( \theta \right) \mbox{d}\mu(\theta)
     }
     {
      \int_{\bR^d} \Big(1 - \cos \langle \xi, x\rangle \Big)\mbox{d}\nu(x)+ 
      \frac{1}{2}\int_{\uS} | \langle \xi, \theta\rangle|^2 \mbox{d}\mu (\theta)
     }.
\end{equation}
It is proved in \cite{BBB} and \cite{BB} that \eqref{symmetric} implies 
\begin{equation}
\|T_mf\|_p\leq (p^*-1)\|f\|_p,\quad 1<p<\infty. 
\end{equation}
This inequality is sharp as these multipliers include $R_2^2-R_1^2$ and $2R_1R_2$. Using Theorem \ref{Non-Symm} we obtain the following related result. 

\begin{theorem}\label{Levyupper} Let $\nu$, $\mu$ be as above and suppose that $\varphi,\,\psi$ take values in $[b, B]$ for some $-\infty<b<B<\infty$.  Then the operator $T_m$ with the symbol \eqref{defmuL}
satisfies
\begin{equation}
\|T_mf\|_p\leq C_{p,b,B}\|f\|_p, \qquad 1<p<\infty.
\end{equation}
\end{theorem} 

Putting $\mu=0$ and using the L\'evy measure $\nu$ of a non-zero symmetric $\alpha$-stable L\'evy process in $\Rd$, with $\alpha \in (0,2)$ (see \cite{BBB} and \cite{BB}), one obtains the multiplier with the symbol
\begin{equation}
 m(\xi)=\frac{\int_{\uS} |\langle \xi,\theta\rangle|^\alpha\phi(\theta)\sigma(\mbox{d}\theta)}{\int_{\uS} |\langle \xi,\theta\rangle|^\alpha\sigma(\mbox{d}\theta)},
\end{equation}
where the so-called spectral measure $\sigma$ is finite and non-zero on $\uS$. By the appropriate choice of $\sigma$ and the use of Theorems \ref{main} and \ref{Levyupper}, we get the following for Marcinkiewicz-type multipliers (see \cite[pp.~109-110]{Ste}).

\begin{corollary}  Let $0<\alpha<2$, $d\geq 2$ and recall that $c_p$ is the Choi constant in \eqref{choi1}. 

(i) For any $J\subsetneq \{1,2,\ldots,d\}$, set
\begin{equation} 
m_{J,\alpha}(\xi)=\frac{\sum_{j\in J} |\xi_j|^\alpha}{\sum_{j=1}^d |\xi_j|^\alpha}.
\end{equation}
Then for $1<p<\infty$,
\begin{equation}
|| T_{m_{J,\alpha}}||_p=C_{p,0,1}=c_p. 
\end{equation}

(ii) Suppose that $d$ is even: $d=2n$, and set
\begin{equation} 
m(\xi)=\frac{|\xi_1^2+\xi_2^2+\ldots+\xi_n^2|^{\alpha/2}}{|\xi_1^2+\xi_2^2+
\ldots+\xi_n^2|^{\alpha/2}+|\xi_{n+1}^2+\xi_{n+2}^2+\ldots+\xi_{2n}^2|^{\alpha/2}}.
\end{equation}
Then for $1<p<\infty$,
\begin{equation}
 ||T_m||_p=C_{p,0,1}=c_p.
 \end{equation}
\end{corollary}

%\begin{remark}
Theorem \ref{main} also gives the lower bound for the norms of the Marcinkiewicz multipliers
$$ m(\xi)=\frac{|\xi_1|^{\alpha_1}|\xi_2|^{\alpha_2}\ldots|\xi_d|^{\alpha_d}}{|\xi|^\alpha},$$
where $\alpha_1,\,\alpha_2,\,\ldots,\,\alpha_d$ are positive numbers and $\alpha=\alpha_1+\alpha_2+\ldots+\alpha_d$, treated in \cite[pp.~109-110]{Ste}. Namely, we have $||T_m||_p\geq C_{p,0,1}=c_p$, for $1<p<\infty$. On the other hand, we have not been able to obtain the reverse bound.

It is also interesting to note here that if $J\subsetneq \{1,2,\ldots,d\}$ and 
\begin{equation} m_J^{log}(\xi)=\frac{\sum_{j\in J}\mbox{ln}(1+\xi_j^{-2})}{\sum_{j=1}^d \mbox{ln}(1+\xi_j^{-2})},
\end{equation}
then 
\begin{equation}
 ||T_{m_J^{log}}||_p\leq C_{p,0,1}=c_p.
 \end{equation}
 Unfortunately these ``logarithmic" multipliers, which arise naturally from the so called tempered stable L\'evy processes (see \cite{BBB}), are not homogeneous of order $0$ and hence the opposite inequality, while still could hold, does not follow from Theorem \ref{main}. 
%\end{remark}

We organize the rest of the paper  as follows. In \S2 we give the proof of the lower $L^p$ bound for multipliers, Theorem \ref{main}.  This proof is a modification of the arguments used by Geiss, Montgomery-Smith and Saksman in  \cite{GMSS}. \S3 is devoted to the proof of Theorem \ref{Non-Symm}: we show there how to deduce \eqref{maininmart} from the discrete martingale inequality \eqref{martin}.  Finally, in \S4 we sketch the proof of the upper bound of Theorem \ref{Riesz1} using the now well known arguments from \cite{BMH}. 

\section{proof of Theorem \ref{main}}
With no loss of generality, we may assume that we have $m(1,0,0,\ldots,0)=b$ and $m(0,1,0,\ldots,0)=B$, rotating and reflecting the multiplier if the equalities do not hold. For the sake of convenience and clarity, we split the proof into several steps.
%\smallskip
\bigskip

\emph{Step 1. The passage from $\R^d$ to the torus $\mathbb{T}^d=(-\pi,\pi]^d$.} Given a smooth and homogeneous multiplier $m$ on $\R^d\setminus\{0\}$, denote by $\tilde{m}$ the corresponding multiplier acting on functions given on $\mathbb{T}^d$. That is, let
\begin{equation}\label{defm}
 T_{\tilde{m}}f(\theta)=\sum_{k\in \mathbb{Z}^d} \hat{f}(k)e^{i\langle k, \theta\rangle}m(k),\quad \theta \in \mathbb{T}^d,
\end{equation}
where, as usual, $\hat{f}(k)=(2\pi)^{-d}\int_{\mathbb{T}^d} e^{-i \langle k,\theta\rangle}f(\theta)\mbox{d}\theta$ and $m(0)=\omega_{d-1}^{-1}\int_{S^{d-1}}m(x)\mbox{d}x$ is the average over the unit sphere in $\R^d$. 

A remarkable fact is that for $1<p<\infty$, the $L^p$ norms of the multipliers $m$ and $\tilde{m}$ coincide. We have the following result due to de Leeuw \cite{L}.
\begin{theorem}
For any $m$ as above and any $1<p<\infty$,
\begin{equation}\label{passage}
 ||T_m:L^p(\R^d)\to L^p(\R^d)||= ||T_{\tilde{m}}:L^p(\mathbb{T}^d)\to L^p(\mathbb{T}^d)||.
\end{equation}
\end{theorem}
Thus it suffices to establish the appropriate lower bound for the norm on the right.
\bigskip
%\smallskip
 
\emph{Step 2. Picking a dyadic martingale and its transform}. Let $f =(f_n)_{n=1}^N$ be a finite, real-valued Paley-Walsh martingale. That is, for 
$n=1,\,2,\,\ldots,\,N$, we have
$$df_n = \e_nd_n(\e_1,\e_2,\ldots,\e_{n-1}),$$
where $\e_1$, $\e_2$, $\ldots$, $\e_N$ is a sequence of independent Rademacher random variables, $d_n:\{-1,1\}^{n-1}\to \R$ are fixed functions, $n=2,\,3,\,\ldots,\,N$, and $d_1$ is a constant. 
Suppose that $\alpha=(\alpha_k)_{k=1}^N$ is a deterministic sequence with each term taking values in $\{b,B\}$ and let $g=(g_n)_{n=1}^N$ be the transform of $f$ by $\alpha$.

%\smallskip
\bigskip

\emph{Step 3. Representing $f$ and $g$ as functions on $(\mathbb{T}^d)^N$.} 
Consider two functions $a^-$, $a^+$ on $\mathbb{T}^d$, defined by $a^-(\theta)=\,$sgn$\,\theta_1$ and $a^+(\theta)=\,$sgn$\,\theta_2$. It is not difficult to see that $T_{\tilde{m}} a^-=ba^-$ and $T_{\tilde{m}}a^+=Ba^+$. Indeed, we easily check that $\widehat{a^-}(k)=0$ if $k_1=0$ or $k_j\neq 0$ for some $j>1$. Consequently, by \eqref{defm}, 
\begin{equation*}
\begin{split}
 T_{\tilde{m}}a^-(\theta)&=\sum_{k_1\in \mathbb{Z}\setminus\{0\}} \widehat{a^-}((k_1,0,0,\ldots,0))e^{ik_1\overline{\theta_1}}m((k_1,0,0,\ldots,0))\\
 &=m(1,0,0,\ldots)\sum_{k_1\in \mathbb{Z}\setminus\{0\}} \widehat{a^-}((k_1,0,0,\ldots,0))e^{ik_1\overline{\theta_1}}\\
 &=b\sum_{k\in\mathbb{Z}}\widehat{a^-}(k)e^{i\langle k,\theta\rangle}=b a^-(\theta).
 \end{split}
 \end{equation*}
The equality $T_{\tilde{m}}a^+=Ba^+$ is proved in the same manner. Now, introduce the sequence $\psi=(\psi_k)_{k=1}^N$ of functions on $\mathbb{T}^d$ by
$$ \psi_k=\begin{cases}
a^- & \mbox{if }\alpha_k=b,\\
a^+ & \mbox{if }\alpha_k=B,
\end{cases}$$
so that
\begin{equation}\label{action}
 T_{\tilde{m}}\psi_k=\alpha_k\psi_k\qquad \mbox{for }k=1,\,2,\,\ldots,\,N.
\end{equation}
We have that $(\psi_1(\theta^1),\psi_2(\theta^2),\ldots,\psi_N(\theta^N))$ has the same distribution (as a function of 
$(\theta^1,\theta^2,\ldots,\theta^N)\in (\mathbb{T}^d)^N$ with normalized measure) as $(\e_1,\e_2,\ldots,\e_N)$. Therefore,
$$ \left(\sum_{k=1}^n \psi_k(\theta^k)d_k\big(\psi_1(\theta^1),\psi_2(\theta^2),\ldots,\psi_{k-1}(\theta^{k-1})\big)\right)_{n=1}^N$$
has the same distribution as the initial martingale $f$. Furthermore, the transform
$g$ can be represented in the form
$$ \left(\sum_{k=1}^n [T_{\tilde{m}}\psi_k](\theta^k)d_k\big(\psi_1(\theta^1),\psi_2(\theta^2),\ldots,\psi_{k-1}(\theta^{k-1})\big)\right)_{n=1}^N,$$
in virtue of \eqref{action}.
\bigskip

%\smallskip

\emph{Step 4. Applying the result of Geiss, Montgomery-Smith and Saksman.} We shall need the following fact. A stronger, Banach-space-valued version appears as Lemma 3.3 in \cite{GMSS}.
\begin{theorem}
Let $1<p<\infty$ and assume that the multiplier $m$ is real and even. For $k\geq 1$, let $E_k$ be the closure in $L^p((\mathbb{T}^d)^k)$ of the finite real trigonometric polynomials
% on $(\mathbb{T}^d)^k$ of the form
$$ \Phi_k(\theta^1,\ldots,\theta^k)
=\sum_{\ell^1\in \mathbb{Z}^d} \ldots \sum_{\ell^k\in \mathbb{Z}^d} e^{i\langle \ell^1,\theta^1\rangle}\ldots 
e^{i\langle \ell^k,\theta^k\rangle} c_{\ell^1,\ldots,\ell^k},$$
such that %Im$(\Phi_k)\equiv 0$, only finitely many of $c_{\ell^1,\ldots,\ell^k}\in \mathbb{C}$ are non-zero, and 
$c_{\ell^1,\ldots,\ell^k}=0$, whenever $\ell^k=0$ (so that $\int_{\mathbb{T}^d}\Phi_k(\theta^1,\ldots,\theta^k)\mbox{d}\theta^k=0$). Let $T_{\tilde{m}}^k$ be an operator on $E_k$, %L^p((\mathbb{T}^d)^k)$
defined on the above polynomials by
$$ (T_{\tilde{m}}^k\Phi_k)(\theta^1,\ldots,\theta^k)
=\sum_{\ell^1\in\mathbb{Z}^d}\ldots\sum_{\ell^k\in
\mathbb{Z}^d} m(\ell^k)e^{i\langle \ell^1,\theta^1\rangle}\ldots 
e^{i\langle \ell^k,\theta^k\rangle}c_{\ell^1,\ldots, \ell^k},$$
for all $\theta^1,\ldots,\theta^k\in \mathbb{T}^d$. Then one has
\begin{equation*}
\begin{split}
\Bigg|\!\Bigg|\sum_{k=1}^N [T_{\tilde{m}}^k\Phi_k](\theta^1,\ldots,\theta^k)&\Bigg|\!\Bigg|_{L^p((\mathbb{T}^d)^N)}\\
&
\leq |\!|T_{\tilde{m}}:L^p(\mathbb{T}^d)\to L^p(\mathbb{T}^d)|\!|
\left|\!\left|\sum_{k=1}^N\Phi_k(\theta^1,\ldots,\theta^k)\right|\!\right|
_{L^p((\mathbb{T}^d)^N)}
\end{split}
\end{equation*}
for all $\Phi_1\in E_1$, $\ldots$, $\Phi_N\in E_N$.
\end{theorem}

Let us apply this result to the representations of $f$ and $g$, setting
$$ \Phi_k(\theta^1,\ldots,\theta^k)=\psi_k(\theta^k)d_k\big(\psi_1(\theta^1),\psi_2(\theta^2),\ldots,\psi_{k-1}(\theta^{k-1})\big)$$
for all $k=1,\,2,\,\ldots,\,N$ and all $\theta^1,\theta^2,\,\ldots,\theta^N\in\mathbb{T}^d$. Then $\Phi_k\in E_k$ for all $k$: the equality $\int_{\mathbb{T}^d} \Phi_k\, \mbox{d}\theta^k=0$ is guaranteed by the martingale property. We obtain
$$ ||g_N||_p \leq ||T_{\tilde{m}}:L^p(\mathbb{T}^d)\to L^p(\mathbb{T}^d)||\,||f_N||_p.$$
Since $N$, $f$ and the transforming sequence $\alpha$ were arbitrary, we get, by \eqref{passage},
$$ ||T_m:L^p(\R^d)\to L^p(\R^d)||= ||T_{\tilde{m}}:L^p(\mathbb{T}^d)\to L^p(\mathbb{T}^d)||\geq C_{p,b,B}.$$
This completes the proof.

\section{Proof of Theorem \ref{Non-Symm}}
First let us first check that the non-symmetric version \eqref{nonsymm} of differential subordination generalizes the martingale transforms by a predictable sequences taking values in $[b,B]$. To do this, let $f$ be a discrete-time martingale and assume that $g$ is its transform by an appropriate sequence $v=(v_n)_{n\geq 0}$. Let us treat $f$, $g$ as continuous-time martingales $X$, $Y$  via the identification $X_t=f_{\lfloor t \rfloor}$ and $Y_t=g_{\lfloor t\rfloor}$, $t\geq 0$. Then both sides of \eqref{nonsymm} are zero for non-integer $t$, and
\begin{eqnarray*}
&&\mbox{d}\left[Y-\frac{b+B}{2}X,Y-\frac{b+B}{2}X\right]_n-\mbox{d}\left[\frac{B-b}{2}X,\frac{B-b}{2}X\right]_n\\\\
&&=dg_n^2-(b+B)df_ndg_n+\frac{(b+B)^2}{4}df_n^2-\frac{(B-b)^2}{4}df_n^2\\\\
&&=(v_n-B)(v_n-b)df_n^2,
\end{eqnarray*}
which is nonpositive when $v_n\in [b,B]$.  Thus \eqref{nonsymm} is satisfied and, in particular, the sharpness in \eqref{maininmart} follows immediately from the passage to discrete-time martingale transforms. 

 To prove \eqref{maininmart}, fix $1<p<\infty$ and note that we may restrict ourselves to $X\in L^p$, since otherwise there is nothing to prove. Then, by Burkholder's inequality \eqref{bur2}, we also have $Y\in L^p$, because $Y$ is differentially subordinate to $(|b|+|B|)X$. Let $V:\R\times \R\to\R$ be the function given by 
$$ V(x,y)=|y|^p-C_{p,b,B}^p|x|^p.$$ 
For any $x,\,y\in \R$, let $M(x,y)$ denote the class of all simple martingale pairs $(f,g)$ starting from $(x,y)$ such that $dg_n=v_ndf_n$, $n\geq 1$, for some deterministic sequence $v$ with terms in  $\{b,B\}$. Introduce the function $U:\R\times \R\to\R$ by
$$ U(x,y)=\sup\{\E V(f_n,g_n)\},$$
where the supremum is taken over all $n$ and all $(f,g)\in M(x,y)$. Of course, $V\leq U$, since the constant pair $(f,g)\equiv (x,y)$ belongs to $M(x,y)$. Furthermore, 
\begin{equation}\label{init}
\mbox{if $y=wx$ for some $w\in [b,B]$, then }
U(x,y)\leq 0.
\end{equation}
This follows from the definition of $U$ and the fact that for such $x,\,y$, the condition $(f,g)\in M(x,y)$ implies that $g$ is the transform of $f$ by a predictable sequence with values in $[b,B]$. Next, using the splicing argument of Burkholder (see e.g. \cite{Bu}) we see that 
\begin{equation}\label{concavity}
\mbox{$U$ is concave along all lines of slope $b$ or $B$.}
\end{equation}
Furthermore, as we shall prove now,
\begin{equation}\label{convexity}
\mbox{for any fixed $x$, the function $U(x,\cdot)$ is convex.}
\end{equation}
 To show this, take any $\lambda\in (0,1)$, $y^-,\,y^+\in \R$ and let $y=\lambda y^-+(1-\lambda)y^+$. Pick any pair $(f,g)\in M(x,y)$ and observe that $(f,g+(y^--y))\in M(x,y^-)$, $(f,g+(y^+-y))\in M(x,y^+)$. Consequently,
\begin{equation*}
\begin{split}
 \E &V(f_n,g_n)
=\E\big[|g_n|^p-C_{p,b,B}^p|f_n|^p\big]\\
 &=\E \big[ |\lambda(g_n+(y^--y))+(1-\lambda)(g_n+(y^+-y))|^p-C_{p,b,B}^p|f_n|^p\big]\\
 &\leq \lambda \E\big[|g_n+(y^--y)|^p-C_{p,b,B}^p|f_n|^p\big]+(1-\lambda)\E\big[|g_n+(y^+-y)|^p-C_{p,b,B}^p|f_n|^p\big]\\
 &\leq \lambda U(x,y^-)+(1-\lambda) U(x,y^+)
 \end{split}
 \end{equation*}
and it suffices to take supremum over $n$ and  $(f,g)$ to get the convexity of $U(x,\cdot)$. Define now 
$\overline{U},\,\overline{V}:\R^2\to \R$ by
$$ \overline{U}(x,y)=U\left(\frac{2}{B-b}x,\frac{B+b}{B-b}x+y\right)$$
and 
$$ \overline{V}(x,y)=V\left(\frac{2}{B-b}x,\frac{B+b}{B-b}x+y\right).$$
We easily check that \eqref{concavity} means that $\overline{U}$ is concave along all lines of slope $\pm 1$ and that \eqref{convexity} carries over to $\overline{U}$. 
Let $\psi:\R\times \R\to [0,\infty)$ be a $C^\infty$ function, supported on the unit ball of $\R^2$, satisfying $\int_{\R^2}\psi=1$. For any $\delta>0$, define $U^\delta,\,V^\delta:\R^2\to \R$ by the convolutions
$$ U^\delta(x,y)=\int_{\R^2} \overline{U}(x+\delta r,y+\delta s)\psi(r,s)\mbox{d}r\mbox{d}s$$ 
and
$$ V^\delta(x,y)=\int_{\R^2} \overline{V}(x+\delta r,y+\delta s)\psi(r,s)\mbox{d}r\mbox{d}s.$$ 
Since $V\leq U$, we have $\overline{V}\leq \overline{U}$ and hence also $V^\delta\leq U^\delta$. Furthermore, the function $U^\delta$ is of class $C^\infty$ and inherits the concavity and the convexity properties of $\overline{U}$. Therefore, we have that
\begin{equation}\label{boundsW}
 U^\delta_{xx}\pm 2U^\delta_{xy}+U^\delta_{yy}\leq 0 \qquad \mbox{and}\qquad U^\delta_{yy}\geq 0\qquad \mbox{on }\,\, \R^2.
\end{equation}
These estimates imply that for all $x,\,y,\,h,\,k\in \R$ we have
$$ U_{xx}^\delta(x,y)h^2+2U_{xy}^\delta(x,y)hk+U_{yy}^\delta(x,y) k^2\leq \frac{U^\delta_{xx}(x,y)-U^\delta_{yy}(x,y)}{2}(h^2-k^2).$$
To see this, we transform the inequality into
$$ (U_{xx}^\delta(x,y)+U_{yy}^\delta(x,y))\frac{h^2+k^2}{2}+2U_{xy}^\delta(x,y)hk\leq 0,$$
and this bound follows easily from \eqref{boundsW} and the trivial estimate $2|hk|\leq h^2+k^2$. Pick two real martingales $X'$, $Y'$ bounded in $L^p$ such that $Y'$ is differentially subordinate to $X'$. 
Then there is a nondecreasing sequence $(\tau_n)_{n\geq 0}$ of stopping times, which converges to $+\infty$ almost surely and $\tau_n$ depends only on $X'$, $Y'$ and $n$, such that
$$ \E U^\delta(X'_{\tau_n\wedge t},Y'_{\tau_n\wedge t})\leq \E U^\delta(X'_0,Y'_0).$$
We refer the reader to Wang \cite{W} for details. 
Since $V^\delta\leq U^\delta$, we get
$$ \E V^\delta(X'_{\tau_n\wedge t},Y'_{\tau_n\wedge t})\leq \E U^\delta(X'_0,Y'_0).$$
Let $\delta\to 0$ and use Lebesgue's dominated convergence theorem to obtain
$$ \E \overline{V}(X'_{\tau_n\wedge t},Y'_{\tau_n\wedge t})\leq \E \overline{U}(X'_0,Y'_0)$$ 
(we note here that the required majorants are of the form $c[(X')^*+(Y')^*]^p$ and their integrability is guaranteed by Doob's maximal inequality.)
 Apply this bound to the pair
$$ X'=\frac{B-b}{2}X\qquad \mbox{ and }\qquad Y'=Y-\frac{B+b}{2}X,$$
and observe that the differential subordination of $Y'$ to $X'$ is equivalent to \eqref{nonsymm}. As the result, we get
$$ \E V(X_{{\tau_n}\wedge t},Y_{\tau_n\wedge t})\leq \E U(X_0,Y_0).$$
However, $U(X_0,Y_0)\leq 0$: use \eqref{init} and the remark below Theorem \ref{Non-Symm}. Therefore,
$$ \E |Y_{\tau_n\wedge t}|^p\leq C_{p,b,B}^p\E|X_{\tau_n\wedge t}|^p$$
and it suffices to first let $n\to \infty$ and then $t\to \infty$  to obtain the desired bound.
%\end{proof}

\section{The Upper Bound in Theorem \ref{Riesz1}}

The upper bound in Theorem \ref{Riesz1} follows immediately from Theorem \ref{Non-Symm} and the stochastic representation for the Riesz transforms as presented in \cite{BMH}.  We also refer the reader to \cite{Ban}, \S3.4,  for a detailed extension of this argument to a wider collection of operators.  Here we only explain how the subordination condition \eqref{nonsymm} of Theorem \ref{Non-Symm} enters into the picture. 
Let $(W,t)$ be  the space-time Brownian motion in $\R^d\times [0,\infty)$. For any sufficiently regular $f$ on $\R^d$, we represent it as the stochastic integral
$$ f \sim X=\int_0^T \nabla U_f(W_{s},{T-s}) \cdot \mbox{d}W_s,$$
where $U_f$ stands for the heat extension of $f$ to the half--space $\R^d\times [0,\infty)$ and $T$ is a large positive number. For a detailed description of what we mean here by the symbol $``\sim"$, see \cite{BMH} or \cite{Ban}.  Then $\mathcal{S}_{\mathbb{A}}$ can be represented as the conditional expectation of the martingale transform of $X$ by $\mathbb{A}$.  That is, 
$$\mathcal{S}_{\mathbb{A}}f(x)\sim E\left[ \,Y\big|\,
W_T=(x,0)\, \right],$$
where 
$$ Y=\int_0^T \mathbb{A}\nabla U_f(W_{s},T-s) \cdot\mbox{d}W_s.$$
Now, if we set  $\xi=\nabla U_f(W_{t},T-t)$, then
\begin{equation*}
\begin{split}
 \mbox{d}\left[Y-\frac{b+B}{2}X,Y-\frac{b+B}{2}X\right]_t&-\mbox{d}\left[\frac{B-b}{2}X,\frac{B-b}{2}X\right]_t\\
 =&\; \left(\left|\mathbb{A}\xi-\frac{b+B}{2}\xi\right|^2-\left|\frac{B-b}{2}\xi\right|^2\right)\mbox{d}t\\
 =&\; \left(|\mathbb{A}\xi|^2-(b+B)(\mathbb{A}\xi,\xi)+bB|\xi|^2\right)\mbox{d}t\\
 =&\;\big\langle(\mathbb{A}-B\mathbb{I})(\mathbb{A}-b\mathbb{I})\xi,\xi\big\rangle\mbox{d}t,
\end{split}
\end{equation*}
where $\mathbb{I}$ stands for the identity matrix of dimension $d$. However, $\mathbb{A}-B\mathbb{I}$ is nonpositive-definite, $\mathbb{A}-b\mathbb{I}$ is nonnegative-definite and the two matrices commute. Hence their product is nonpositive-definite and hence \eqref{nonsymm} is satisfied. Consequently, by \eqref{maininmart},
$$ ||Y||_p\leq C_{p,b,B}||X||_p,$$
which, by the ``transference method," as explained in \cite{BMH} and \cite{Ban},  yields
$$ ||\mathcal{S}_{\mathbb{A}}f||_p\leq C_{p,b,B}||f||_p$$
and hence $||\mathcal{S}_{\mathbb{A}}:L^p(\R^d)\to L^p(\R^d)||\leq C_{p,b,B}$. 
%\end{proof}

\begin{remark}  It is worth observing here that, as the proof of the upper bound of Theorem \ \ref{Riesz1} shows, if we take a real variable coefficient  $d\times d$ symmetric matrix ${\mathbb{A}}(x, t)$,  $x\in \R^d$, $t>0$, with the property that for all $\xi\in \R^d$, $$b|\xi|^2\leq \langle{\mathbb{A}}(x, t)\xi, \xi\rangle\leq B|\xi|^2,$$
for all $(x, t)$, and define the operator 
$$\mathcal{S}_{\mathbb{A}}f(x)\sim E\left[ \,Y\big|\,
W_T=(x,0)\, \right],$$
where this time 
$$ Y=\int_0^T \mathbb{A}(W_s, T-s)\nabla U_f(W_{s},T-s) \cdot\mbox{d}W_s,$$
we get 
$$ ||\mathcal{S}_{\mathbb{A}}f||_p\leq C_{p,b,B}||f||_p, \qquad 1<p<\infty. $$
For more on these variable coefficient ``projections of martingale transforms," see \cite{Ban} and especially Remark 3.4.2 there.
\end{remark}

\end{document}